\documentclass[]{elsarticle}

\usepackage[utf8]{inputenc}
\usepackage{amsmath,amsfonts,amsthm}
\usepackage[linesnumbered,ruled]{algorithm2e}
\usepackage[english]{my-shortcuts}
\usepackage[pagewise]{lineno}

\journal{Journal of Computational and Applied Mathematics}

\begin{document}

\begin{frontmatter}

\title{A note on computing the Smallest Conic Singular Value}

\author[NPL]{Stephane Chretien\corref{cor1}}
\address[label1]{National Physical Laboratory,
Hampton Road, Teddington, TW11 0LW, UK}

\cortext[cor1]{Corresponding author}

\ead{stephane.chretien@npl.co.uk}
\ead[url]{https://sites.google.com/view/stephanechretien/home}

\begin{abstract}
The goal of this note is to study the smallest conic singular value of a matrix from a Lagrangian duality viewpoint and provide an efficient method for its computation.
\end{abstract}

\begin{keyword}
Conic singular value, subdifferential, Lagrange duality, Newton's method.
\end{keyword}

\end{frontmatter}


	\section{Introduction} 
	Given a cone $K$, the minimum conic singular value of a matrix $A\in \mathbb R^{d\times n}$ is defined as 
	\begin{align*}
	\sigma_{\min}(A;K) & = \min_{x\in K,\: \Vert x\Vert_2=1} \quad \Vert A x\Vert_2.
	\end{align*} 
	The conic singular value has emerged as an important concept in statistics \cite{AmelunxenetAl:InfInf14}, \cite{Tropp:Renaissance15} and \cite{ThrampoulidisHassibi:ISIT15}. The associated problem of computing the smallest conic eigenvalue of a matrix $M$ has also been extensively studied in the field of optimisation theory in relation to $K$ co-positivity \cite{hiriart2010variational}, \cite{da2010cone}, \cite{seeger2003eigenvalues}, etc. Application to contact problems \cite{da2010cone}, elliptic PDE's \cite{riddell1977eigenvalue}, 
	problems on matrix cones \cite{seeger2003eigenvalues} are numerous.
	Some interesting methods have been proposed for the computation of cone-constrained eigenvalue problems, e.g. the SPA method of \cite{da2009numerical}.  
	
	The goal of this short note is to describe an efficient dual method for solving this problem.One problem with the computation of the minimum conic singular value is that the constraint set is nonconvex due to the spherical constraint. In order to circumvent this problem, we use a Lagrangian framework and Lagrangian duality and a well known formula for the solution of the problem of quadratic programming over the sphere. 
	
	The plan  of the paper is as follows. In section \ref{lag}, the Lagragian approach is explained and a strong duality result is proved. An efficient dual algorithm is described in Section \ref{comp}. Simulation experiments showing the efficiency of the dual approach is tested are presented in Section \ref{sim}.

	\section{The Lagrangian approach}
	\label{lag}
	In the present section, we study the smallest singular value from a dual viewpoint, following the approach of \cite{HULLII}.
	
	\subsection{Presentation of the method} 
	
	We will address the minimum singular value problem via the solution of  
	\begin{align}
	\frac12 \ \sigma_{\min}(A;K)^2 & =\min_{x\in K,\: \Vert x\Vert_2=1} \quad \frac12 \Vert A x\Vert_2^2.
	\label{sqsing} 
	\end{align} 
	For this purpose, we form the Lagrangian function \cite[Chapter XII]{HULLII}
	\begin{align*}
	L(x,u) & = \frac12 \Vert A x\Vert_2^2 + \langle u,x\rangle 
	\end{align*}
	and the dual function 
	\begin{align}
	\theta(u) & = \min_{\Vert x\Vert_2=1} \quad L(x,u).
	\label{dualfunc} 
	\end{align}
	The dual problem is 
	\begin{align*}
	\sup_{K\in K^\circ} \quad \theta(u).
	\end{align*}
	In order to motivate the dual approach it is common to point out that 
	\begin{align*}
	\sup_{u \in K^\circ} \quad L(x,u) & = 
	\begin{cases}
	+ \infty \textrm{ if } \quad x \not \in K \\
	\\
	\frac12 \Vert A x\Vert_2^2  \quad \textrm{ otherwise}.
	\end{cases}
	\end{align*}
	Thus, problem \eqref{sqsing} is equivalent to 
	\begin{align}
	\inf_{\Vert x\Vert_2=1} \quad \sup_{u\in K^\circ} \quad L(x,u).
	\label{primal}
	\end{align} 
	On the other hand, the dual problem of maximizing $\theta$ is equivalent to solving 
	\begin{align}
	\sup_{u\in K^\circ} \quad \inf_{\Vert x\Vert_2=1} \quad  L(x,u),
	\label{dual}
	\end{align} 
	with the hope that the two optimal values coincide. In general, weak duality holds, i.e.
	\begin{align}
	    {\rm OPT}\ \textrm{\eqref{primal}} & \ge {\rm OPT}\ \textrm{\eqref{dual}}
	\end{align}
	where "OPT" denotes the optimal value. When equality holds, we say that the problem enjoys the strong duality property. Strong duality automatically holds for linear programs. It
	also holds for convex programs when an additional constraint qualification condition holds
	\cite[Chapter XII]{HULLII}. It may not hold for nonconvex programs except in certain circumstances. For instance, if $K=\mathbb R^n$ in \eqref{sqsing}, then, strong duality 
	holds despite the nonconvex spherical constraint. 
	
	\subsection{The subdifferential of $\theta$} 
	
	The theory of \cite[Chapter XII]{HULLII} directly provides  a complete description of the subdifferential of $\theta$, because the so called {\em filling property} holds due to the compactness of the spherical constraint. More precisely,
	\begin{align*}
	\partial \theta (u) & = \overline{{\rm conv}} \left( \mathcal X^*_u \right)
	\end{align*} 
	where $\mathcal X^*_u$ is the set of minimizers in \eqref{dualfunc} and $\overline{{\rm conv}}$ denotes the closure of the convex hull. 
	
	Lemma 2.2 in \cite{Hager:SIAMOpt01} gives the explicit form of the solution set $\mathcal X^*_u$. Let $\lambda_1 \le \ldots \le\lambda_n$ be the eigenvalues of $A^tA$ and $\phi_1$,\ldots,$\phi_n$ be associated pairwise orthogonal, unit-norm eigenvectors. 
	Let $\gamma_i= \phi_i^tu$, $i=1,\ldots,n$. Let $\mathcal E_1=\{i \text{ s.t. } \lambda_i=\lambda_1 \}$ and $\mathcal E_+=\{i \text{ s.t. } \lambda_i>\lambda_1 \}$. Then, 
	$x^*_u$ belongs to $\mathcal X^*_u$ if and only if 
	\begin{align}
	x^* & = \sum_{i=1}^n c^*_i \phi_i
	\label{hager}
	\end{align}
	and 
	\begin{enumerate}
		\item {\em degenerate case:} if $\gamma_i=0$ for all $i \in \mathcal E_1$ and 
		\begin{align*}
		\sum_{i\in \mathcal E_+} \: \frac{\gamma_i^2}{(\lambda_i-\lambda_1)^2} \le 1. 
		\end{align*} 
		then 
		\begin{align}
		c_i^*=\gamma_i/(\lambda_i-\lambda_1), 
		\label{cdeg}
		\end{align}
		$i\in \mathcal E_+$ and $c_i^*$, $i\in \mathcal E_1$ are arbitrary under the constraint 
		that 
		\begin{align}
		\sum_{i\in \mathcal E_1} \quad c^{*^2}_i = 1-\sum_{i\in \mathcal E_+} \quad c^{*^2}_i. 
		\label{sphere}
		\end{align}
		
		\item {\em nondegenerate case:} if not in the degenerate case, 
		\begin{align}
		    c_i^*=\gamma_i/(\lambda_i-\mu), 
		    \label{cnondeg}
		\end{align}
		$i=1,\ldots,n$ for $\mu > -\lambda_1$ which is a solution of 
		\begin{align}
		    \sum_{i=1,\ldots,n} \: \frac{\gamma_i^2}{(\lambda_i-\mu)^2} = 1. 
		\label{mu} 
	\end{align} 
	\end{enumerate}
	Moreover, lower and upper bounds are given in \cite{Hager:SIAMOpt01} and the value of $\mu$ can be found very quickly by the bisection method.
	
	Using this representation of $\mathcal X^{c^*}_u$, we immediately deduce the following result.
	\begin{prop}
		The set $\mathcal X^{c^*}_u$ is a unit sphere of dimension less than or equal to $n$.
		\label{sphinside} 
	\end{prop}
	
	\subsection{Strong duality for the minimum conic singular value} 
	
	Consider now a full Lagrangian scheme for this problem. For this purpose, let us define
	\begin{align*}
	L^c(x,u,v) & = \frac12 \Vert A x\Vert_2^2 + \langle u,x\rangle + v\ \Big(\Vert x\Vert_2^2-1\Big) 
	\end{align*}
	and the dual function 
	\begin{align}
	\theta^c(u,v) & = \min_{\Vert x\Vert_2\le 1} \quad L^c(x,u,v).
	\label{fulldualfunc} 
	\end{align}
	
	One interesting point with this full Lagrangian approach is that we can easily prove that strong duality holds. 
	\begin{theo}
		Strong duality holds for the dual problem 
		\begin{align*}
		\sup_{(u,v) \in K^\circ \times \mathbb R} \quad \theta^c(u,v) & = \sigma_{\min}(A;K)^2.
		\end{align*}
	\end{theo}
	\begin{proof}
		Indeed the subdifferential of $\theta^c$ can be computed easily as before as
		\begin{align*}
		\partial \theta^c (u) & = \overline{{\rm conv}} \left(
		\left[
		\begin{array}{c}
		\mathcal X^{c^*}_u \\
		\\
		\Vert \mathcal X^{c^*}_u \Vert_2^2-1
		\end{array}  
		\right] \right)
		\end{align*} 
		where $\mathcal X^{c^*}_u$ is the set of minimizers in \eqref{fulldualfunc}. 	
		One interesting feature of this subdifferential is that it has one (and only one) quadratic component. From this, we can deduce that the image of the map 
		\begin{align*}
		x & \mapsto 
		\left[
		\begin{array}{c}
		x \\
		\\
		\Vert x \Vert_2^2-1
		\end{array}  
		\right]
		\end{align*}
		On the other hand, the optimality condition for the dual problem is that there exists $g$ in the normal cone to $K^\circ$ at $u^*$ such that 
		\begin{align*}
		(g,0) & \in \partial \theta^c(u^*).
		\end{align*} 
		In other words, 
		\begin{align*}
        \left[\begin{array}{c}
	    g \\
	    0
	    \end{array}
	    \right]
        & \in \overline{{\rm conv}} \left(
		\left[
		\begin{array}{c}
		\mathcal X^{c^*}_{u^*} \\
		\\
		\Vert \mathcal X^{c^*}_{u^*} \Vert_2^2-1
		\end{array}  
		\right] \right).
		\end{align*} 
		Since $\mathcal X^{c^*}_{u^*}$ is included in the unit ball by the definition of $\theta^c$, the set 
		\begin{align*}
		\mathcal S_{u^*} & = \left[
		\begin{array}{c}
		\mathcal X^{c^*}_{u^*} \\
		\\
		\Vert \mathcal X^{c^*}_{u^*} \Vert_2^2-1
		\end{array}  
		\right]
		\end{align*} 
		is included in a compact set and therefore. Thus, the closure of its convex hull is the convex hull of its closure. Because the Lagrangian is continuous, $\mathcal X^{c^*}_{u^*}$
		is also compact, the set $\mathcal S_{u^*}$ is equal to its closure. Now since $\mathcal X^{c^*}_{u^*}$ 
		is a sphere by Corollary \ref{sphinside}, Brickman's celebrated theorem \cite{Barvinok:AMS02} 
		on the quadratic image of a sphere gives that $\mathcal S_{u^*}$ is convex. Therefore, the optimality 
		condition becomes that there exists $g$ in the 
		normal cone to $K^\circ$ at $u^*$ such that 
		\begin{align*}
    		(g,0) & \in \mathcal S_{u^*}
		\end{align*}
		which implies that there exists $x^*_{u^*}$ in $\mathcal X^{c^*}_{u^*}$ such that  
		\begin{align*}
		\left[
		\begin{array}{c}
	        g \\
	        0
	    \end{array}
	    \right]
         & =
		\left[
		\begin{array}{c}
		x^*_{u^*} \\
		\\
		\Vert x^*_{u^*} \Vert_2^2-1
		\end{array}  
		\right].
		\end{align*}  
		The end of the proof is standard. Since $g$ is in the normal cone to $K^\circ$ at $u^*$, and $g=x_{u^*}$, we have 
		\begin{align*}
		\langle u^*,x^*_{u^*}\rangle & =0. 
		\end{align*}  
		Thus, 
		\begin{align*}
		\theta^c(u^*) & = \frac12 \Vert Ax_{u^*} \Vert_2^2.
		\end{align*}
		This proves that 
		\begin{align*}
		\theta(u^*) & \ge \frac12 \sigma_{\min} (A;K)^2
		\end{align*}
		but since, by weak duality, 
		\begin{align*}
		\theta(u^*) & \le \frac12 \sigma_{\min} (A;K)^2
		\end{align*}
		we finally obtain that strong duality holds and that $x_{u^*}$ is a solution. 
	\end{proof}
	We now obtain the following corollary. 
	\begin{cor}
		Strong duality holds for the dual problem 
		\begin{align*}
		\sup_{(u) \in K^\circ } \quad \theta(u) & = \frac12 \ \sigma_{\min}(A;K)^2.
		\end{align*}
	\end{cor}
	\begin{proof}
		Use the same proof as in \cite{LemarechalOustry:AdvancesInConvexAnalysis01} to obtain that the full Lagrangian scheme 
		is weaker than the Lagrangian scheme. Since the full Lagrangian scheme is exact, so is the Lagrangian one.  
	\end{proof}
	
	\section{Computation of the smallest conic singular value}
	\label{comp}
	\subsection{The polar cone $K^\circ$ of a cone $K$}
	\subsubsection{The polyhedral case}
	The polar cone $K^\circ$ to $K$ is a key object in our computations. It has been studied in \cite{dobler1994matrix}. A simplified version of the main theorem of \cite{dobler1994matrix} is the following.
	
	\begin{theo}\cite[Theorem 4.2]{dobler1994matrix}
	    Assume that the cone $K$ has representation 
	    \begin{align}
	        K & = \left\{ x \in \mathbb R^{n} \mid A^t x \le 0 \textrm{ and } B^t x = 0\textrm{ for } A \in \mathbb R^{n \times m} \textrm{ and  } B \in \mathbb R^{r\times n}\right\}.
	    \end{align}
	    Assume that there exist two matrices $G_A$ and $G_B$ such that $G_AA=I$ and $G_BB=I$. If in addition, $G_AB=0$, 
	    then, we have 
	    \begin{align}
	        K^\circ & = \left\{ y \in \mathbb R^{n} \mid -G_Ay \le 0  \textrm{ and  }(I-AG_A-BG_B)y=0\right\}. 
	    \end{align}
	\end{theo}
	
	Therefore, in order to be able to use this theorem, we have to solve the system of equations $G_AA=I$, $G_BB=I$ and $G_AB=0$. As noticed in \cite{dobler1994matrix}, this system may have more than one solution and these solution can often by found in the form $G_A=(A^tDA)^{-1}A^tD$ for some full rank matrix $D$. If $B=0$, then a solution may be found in the simple form $G_A=(A^tA)^{-1}A^t$.  
	
	\subsubsection{The non-polyhedral case}
	In the nonpolyhedral case, the computation of the polar may be more difficult except for some standard cases where the polar is well known, like for the Positive Semi-Definite cone. One option to get around this is to find a good polyhedral approximation and apply the formulas from the previous section. We will not enter the details here, but will refer instead the interested reader to \cite{bronshteyn1975approximation}.
	
	\subsection{A quasi-Newton approach for optimising the dual problem}
	The dual approach assumes that we can easily compute the polar to the cone $K$. This can be done efficiently using 
	the results in \cite{dobler1994matrix} in the case of polyhedral cones.
	One very efficient method for this type of nonsmooth problem is the bundle method \cite{HULLII}. We will not describe this method here and refer the reader to \cite{HULLII} for an very pedagogical introduction. 
	In Algorithm \ref{dualcomp}, we propose a very simple non-smooth quasi-Newton approach inspired from \cite{yu2010quasi}. 
	
	Before presenting the algorithm, let us recall that the BFGS quasi-Newton update is given by 
	\begin{align}
	    H^{(l+1)} & = H^{(l)} +\frac{yy^t}{y^ts}
	    - \frac{H^{(l)}ss^t H^{(l)}}{s^tH^{(l)}s}
	    \label{qNupdate}
	\end{align}
	with $u^{(l+1)}-u^{(l)}$ and $y=\nabla \theta(u^{(l+1)})-\nabla \theta(u^{(l)})$. This gives the inverse formula: 
	\begin{align}
	    H^{(l+1)^{-1}} & =\Big(I-\frac{sy^t}{y^ts}\Big) H^{(l)^{-1}}\Big(I-\frac{ys^t}{y^ts}\Big) +\frac{ss^t}{y^ts}.
	    \label{qNupdateinv}
	\end{align}

	\vspace{.5cm}

\begin{algorithm}[H]
\SetAlgoLined

Initialise with $u^{(1)}=u^{(0)}=0$ and $H^{(1)}=I$.
 
 \textbf{While $\Vert u^{(l+1)}-u^{(l)} \Vert_2\ge \epsilon$}

 \quad - Compute $\mathcal X^*_{u^{(l)}}$ using 	\eqref{hager} with \eqref{cdeg} and \eqref{sphere} in the degenerate case and 
	        \eqref{cnondeg} and \eqref{mu} in the nondegenerate case.
	        
	         \quad -  Choose $x^*_{u^{(l)}} \in \mathcal X^*_u$ such that 
			\begin{align}
			    H^{(l)^{-1}}\ g^{(l)} & = H^{(l)^{-1}}\ x^*_{u^{(l)}}
			\end{align} 
			is a descent direction for $\theta$.
	        
	         \quad -  Solve 
			\begin{align}
			u^{(l+1)} & \in \textrm{argmin}_{u \in K^\circ}  \ \langle g^{(l)},u\rangle+\frac12 \ \langle u,H^{(l)^{-1}}u\rangle.
			\end{align}
			
			 \quad -  Update $H^{(l+1)^{-1}}$ using the quasi-Newton update rule \eqref{qNupdateinv}
 
\textbf{End while}
 
\textbf{Output}: $x^{(L)}=x^*_{u^{(L)}}$.

\caption{The dual quasi-Newton method \label{dualcomp}}
\end{algorithm}


	\subsection{Simulation experiments}
	\label{sim}
	We ran our method on random instances in dimensions 15000, 20000 and 25000. In these experiments, the conic constraints were sampled randomly from 100 random inequalities with i.i.d. standard Gaussian distribution and the matrix $A$ was chosen randomly from an i.i.d. standard Gaussian distribution. The stopping criterion was set to $\epsilon = 1e-4$.
	
	\vspace{.3cm}
	
	The following table shows the average computational time for different problem sizes both for our method and the SPA method of \cite{da2009numerical}. The symbol '*' stands for 'did not converge in less than 10 times the computation time of our method'.

\begin{center}
\begin{tabular}{c|c|c|c}
\hline
    Problem size & 15000  & 20000 & 25000 \\
\hline    
    Ave. Comp. time our method (sec) & 4167 & 9563 & 24356 \\
    Ave. Comp. time SPA (sec) & * & * & * \\ 
    \hline
\end{tabular}
\end{center}

\begin{figure}[!t]
\centering
\includegraphics[width = 12cm]{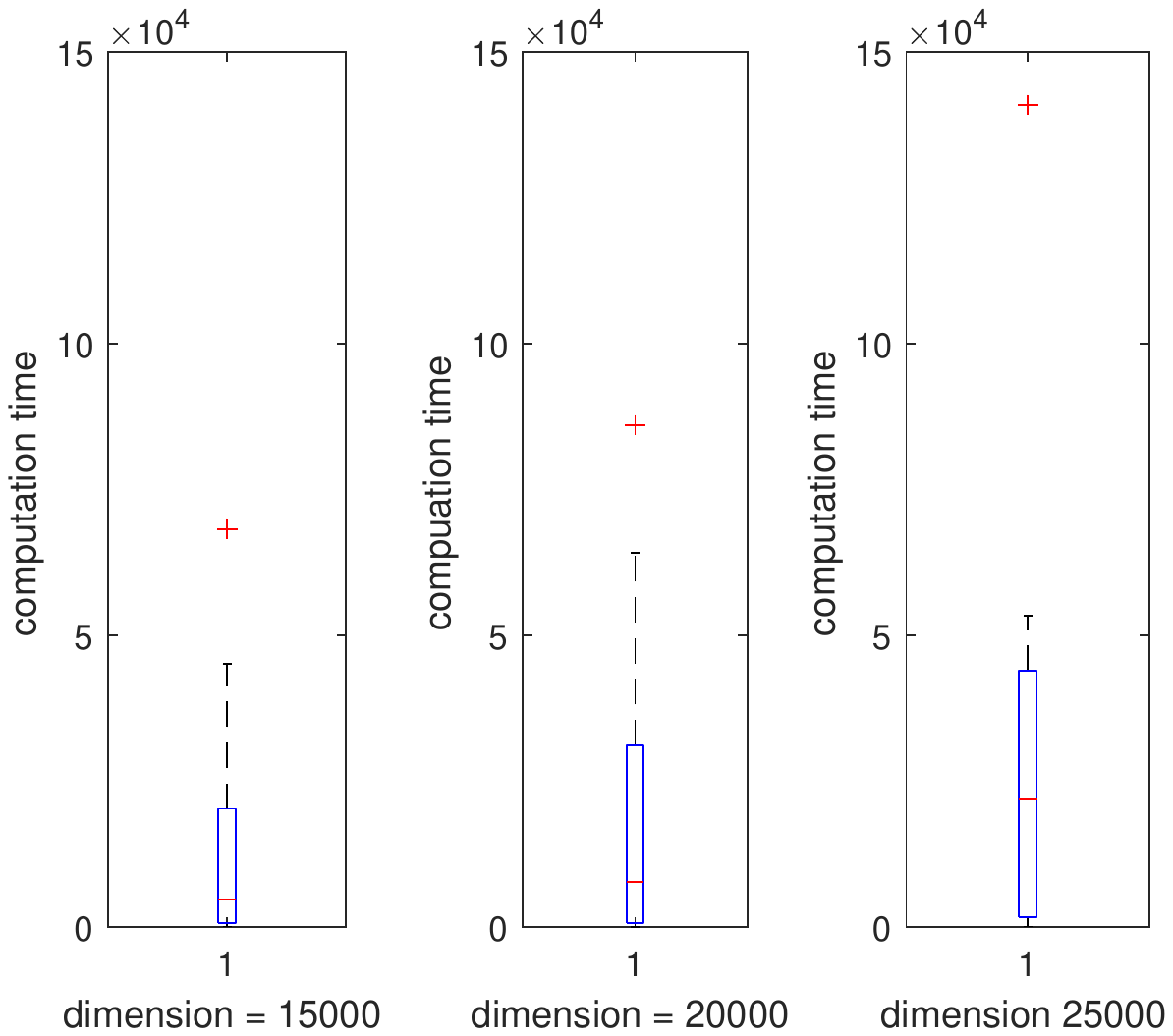}
\caption{
\label{numres} Dispersion of the computation time for 3 different (large) dimensions.} 
\end{figure}

Figure \ref{numres} shows how spread the computation time is around the mean. It also shows that some outliers have a much larger computational time. 

\section{Application to sensor placement in power grids}
In the present section, we describe an application of the conic singular value to the problem of sensor placement in power grids. 

\subsection{Background on the state estimation problem}
Power grids have been a topic of extensive interest lately in the optimisation community. One of the main reason for this surge of activity in the field of power networks stems from the recent interest in estimating the state of the network with as few measurements as possible due to the potential inaccessibility of certain points on the grid. 

From a mathematical viewpoint, our goal is to estimate a vector $V \in \mathbb C^{n}$ whose components are the voltage values at the various buses on the grid. This vector is observed through noisy quadratic measurements $z_l$, $l=1,\ldots,L$ and are related to $V$ via the following equation.
\begin{align}
	z_l & = V^*H_lV+\epsilon_l 
\end{align}
where $\epsilon_1,\ldots,\epsilon_L$ are i.i.d. random variables with distribution $\mathcal N(0,\sigma^2)$. The reason for these measurements to be quadratic stems from the fact that often times, only power measurements are possible due to budget constraints.

\subsection{The least squares problem and a Semi-Definite Relaxation} 
The least squares estimation problem is thus given by 
\begin{align} 
\min_{V\in \mathbb C^{N\times N}} \quad \sum_{l=1}^L \Big(z_l-{\rm trace}(H_lVV^*)\Big)^2
\label{ls}
\end{align} 
Using the change of variable $W=VV^*$, we then have the equivalent Rank constrained Semi-Definite Program 
\begin{align*} 
\min_{W\in \mathbb C^{N\times N}} \quad \sum_{l=1}^L \Big(z_l-{\rm trace}(H_l W)\Big)^2 
\hspace{.5cm} \textrm{ s.t. } W \succeq 0 \hspace{.5cm}  \textrm{ and } \hspace{.5cm}  {\rm rank}(W)=1. 
\end{align*} 
A standard way to obtain a Semi-Definite relaxation is just to relax the rank one constraint. 
The resulting Semi-Definite Relaxation is given by 
\begin{align} 
\min_{W\in \mathbb C^{N\times N}} \quad \sum_{l=1}^L \Big(z_l-{\rm trace}(H_l W)\Big)^2 
\hspace{.5cm}  \textrm{ s.t. } \hspace{.5cm}  W \succeq 0. 
\label{sdr}
\end{align} 

\subsection{A real relaxation}
In Section \ref{pert}, we intend to study a perturbed version of the Semi-Definite Relaxation. Sensitivity analysis of Semi-Definite Programs has been proposed in a number of 
works; see \cite{shapiro1997first}, \cite{bonnans2013perturbation}, \cite{freund2007nonlinear}. However, we are not aware of any sensitivity analysis of complex Semi-Definite Programs. Fortunately, a simple transformation can 
be performed in order to obtain a real Semi-Definite Program from a complex one; 
see \cite{goemans1995improved}. As expected, it suffices to decompose 
the problem into real and complex parts. 

Let $H_l^{\mathcal R}$ (resp. $H_l^{\mathcal I}$) denote the real (resp. imaginary) part of 
$H_l$ $l=1,\ldots, L$. Similarly, let $W^{\mathcal R}$ (resp. $W^{\mathcal I}$) denote the real (resp. imaginary) part of $W$. Let us now define the new variable $\mathcal W$ by 
\begin{align*} 
\mathcal W & = 
\left[
\begin{array}{cc} 
	W^{\mathcal R} & W^{\mathcal I} \\
	-W^{\mathcal I} & W^{\mathcal R}
	\end{array} 
\right].
\end{align*}  
Similarly, let 
\begin{align*} 
\mathcal H_l & = 
\left[
\begin{array}{cc} 
H_l^{\mathcal R} & H_l^{\mathcal I} \\
-H_l^{\mathcal I} & H_l^{\mathcal R}
\end{array} 
\right]
\end{align*}  
$l=1,\ldots,L$. Notice that the matrices $\mathcal W$ and $\mathcal H_l$ enjoy some 
obvious symmetries. Let us denote by $\mathbb M \subset \mathbb R^{2N\times 2N}$ the space of such matrices.  
Finally, for the sake of reducting the redundency, 
let us define the transformation $\mathcal T: \mathbb M 
\mapsto \mathbb R^{N^2+N(N+1)/2}$ as the one of concatenating 
the first $N$ components of the successive columns of the upper triangular part of a given real matrix into one single column vector and let $\mathcal T^{-1}: \mathbb R^{N^2+N(N+1)/2} \mapsto \mathbb M$ denote its inverse. 
Let $w=\mathcal T (\mathcal W)$ and let $h_l=\mathcal T(\mathcal H_l)$, $l=1,\ldots,L$. 
\begin{lem}
Problem (\ref{sdr}) can be transformed into 
\begin{align}
\label{quadsdp} 
\min_{w \in \mathbb R^{2N(2N+1)/2}} \quad \frac12  \left\Vert q- Q w \right\Vert_2^2 \quad \textrm{ s.t. } \quad \mathcal T^{-1}(w) \succeq 0 
\end{align}  
with 
\begin{align*} 
Q = 2 \sum_{l=1,\ldots,L} (\mathcal T^{-1})^* (\mathcal H_l) (\mathcal T^{-1})^* (\mathcal H_l)^t  & \hspace{1cm} \textrm{ and } &  q = 2 \sum_{l=1,\ldots,L} z_l(\mathcal T^{-1})^* (\mathcal H_l)
\end{align*}
\label{transfTau}
\end{lem}
\begin{proof}
Straightforward computation.
\end{proof}
If this Semi-Definite Program has a solution $w^{opt}$ such that $\mathcal T^{-1}(w^{opt})$ has rank one, then,
\begin{align*} 
\mathcal T^{-1}(w^{opt}) & = v^{opt} v^{opt^t}
\end{align*}  
and,  clearly,  $v^{opt}$ is a solution to (\ref{ls}).

\subsection{Sensitivity of the solution} 
\label{pert}
Sensitivity of Nonlinear Semi-Definite Programs such as (\ref{quadsdp}) have been studied in many 
previous works; see e.g. \cite{freund2007nonlinear} and 
the comprehensive \cite{bonnans2013perturbation}. We can also use the main result from \cite{self1987asymptotic} in the particular setting of our Semi-Definite Relaxation.

\begin{theo} \textbf{\cite[(Corollary of) Theorem 2]{self1987asymptotic}.}
	Let 
	\begin{align*}
	l_L(w) & = \sum_{l=1}^L \Big(z_l- \left\langle \left(\mathcal T^{-1}\right)^*)(\mathcal H_l), w \right\rangle \Big)^2
	\end{align*}
	and
	\begin{align*}
	\mathcal C & = \left\{w \mid \mathcal T(w) \succeq 0 \right\},	
	\end{align*}
	and let ${\rm Tang}_\mathcal C(w)$ denote the tangent cone to $\mathcal C$ at $w$
	Consider the optimization problem
	\begin{align} 
	w^*_L & = {\rm argmin}_{w \in \mathbb R^{2N(2N+1)/2}} \quad l_L(w)
	\hspace{.5cm}  \textrm{ s.t. } \hspace{.5cm}  w \in \mathcal C. 
	\label{sdr2}
	\end{align}
	Assume that for all $l=1,\ldots,L$, 
	\begin{align*}
	z_l & =  \left\langle \left(\mathcal T^{-1}\right)^*)(\mathcal H_l),w_0\right\rangle+n_l
	\end{align*} 
	for some vector $w_0 \in \mathbb R^{2N(2N+1)/2}$ satisfying the constraint $\mathcal T(w_0) \succeq 0$ and with $(n_l)_{l=1,\ldots,L}$ an i.i.d. sequence of random variables with distribution $\mathcal N(0,\sigma^2)$. 
	Let 
	\begin{align*}
	I(w) & = \nabla^2 l_L(w).
	\end{align*} 
	Let $G$ be a Gaussian $\mathcal N(0,I(w_0)^{-1})$ vector and let $W$ be the (random) solution of the following optimization problem
	\begin{align*}
	\min_{\tau \in {\rm Tang}_\mathcal C(w_0)} \ \frac12 (Z-\tau) I(w_0) (Z-\tau).
	\end{align*}
	Then,
	\begin{align*}
		\sqrt{L} \ \left(w_L^* - w_0\right) \rightarrow^{\mathcal L} \ W.
	\end{align*}
\end{theo}

\subsection{Optimal design via conic eigenvalue maximisation}

Now that the estimation problem has been linearised via Semi-Definite Relaxation and transformed into a real Semi-Definite Program, we can address the problem of finding an optimal design of experiments. 
For this purpose, we define a selection vector $\delta \in \{0,1\}^n$ whose components will specify if a measurement $z_l$ is being made (if $\delta_l=1$) or not (if $\delta_l=0$). 
Let 
\begin{align*}
l_{\delta,L}(w) & = \sum_{l=1}^L \Big(\delta_l z_l- \left\langle \left(\mathcal T^{-1}\right)^*)(\delta_l \mathcal H_l), w \right\rangle \Big)^2
\end{align*}
and, as in the previous section,
\begin{align*}
\mathcal C & = \left\{w \mid \mathcal T(w) \succeq 0 \right\}.	
\end{align*}
With these notations, the estimation problem becomes the one of solving
\begin{align} 
w^*_{\delta,L} & = {\rm argmin}_{w \in \mathbb R^{2N(2N+1)/2}} \quad l_{\delta,L}(w)
\hspace{.5cm}  \textrm{ s.t. } \hspace{.5cm}  w \in \mathcal C. 
\label{sdr4}
\end{align}

Let 
\begin{align*}
I_{\delta,L}(w) & = \nabla^2 l_{\delta,L}(w).
\end{align*} 
Since, due to the fact that $l_{\delta,L}(w)$ is quadratic, $I_{\delta,L}(w)$ does not depend on $w$, we will simply denote it by 
$I_{\delta,L}$.
A straightforward computation (taking into account that $\delta_l^2=\delta_l$) gives 
\begin{align*}
I_{\delta,L} & = 2 \sum_{l=1}^L \delta_l \ \left(\mathcal T^{-1}\right)^*)(\mathcal H_l)\left(\mathcal T^{-1}\right)^*)(\mathcal H_l)^t.
\end{align*}

Recalling that the asymptotic distribution of $w^*_{\delta,L}$ depends on the behavior of the quadratic
form $\tau^t I_{\delta,L} \tau$ when $\tau$ is subject to lie in the cone ${\rm Tang}_{\mathcal C}(w_0)$,
one of the most adequate objective functions to use is 
the function $\sigma_{\min}(I_{\delta,L},{\rm Tang}_{\mathcal C})$ which denotes the smallest 'conic' singular value of $I_{\delta,L}$, i.e. 
\begin{align}
\min_{\stackrel{\tau \in {\rm Tang}_{\mathcal C}(w_0)}{\Vert \tau \Vert_2=1}} \ \tau^t I_{\delta,L} \tau.
\label{conesingI}
\end{align}
Furthermore, one can easily prove that 
\begin{align*}
{\rm Tang}_{\mathcal C}(w) & =  
\big\{ \tau \mid  \mathcal T^*\left(P_{\mathbb M} (u_ju_j^t)\right)^t\tau \ge 0, \textrm{ for all } j   \\
& \hspace{1cm} \textrm{ for } (u_j)_{j=1,\ldots,{\rm dim}({\rm Ker}(\mathcal T(w))) } \textrm{ a basis of } {\rm Ker}(\mathcal T(w))\big\}. 
\end{align*}
Given these computations, we obtain that 
\begin{align*}
\sigma_{\min}(M,{\rm Tang}_{\mathcal C})
& = \min_{\Vert \tau \Vert_2=1} \ \tau^t M \tau
\end{align*}
under the additional constraint 
\begin{align*}
\mathcal T^*\left(P_{\mathbb M} (u_ju_j^t)\right)^t\tau \ge 0, 
\end{align*}
for all $j$ and  
$(u_j)_{j=1,\ldots,{\rm dim}({\rm Ker}(\mathcal T(w))) } \textrm{ a basis of } {\rm Ker}(\mathcal T(w))$.

The optimal choice of $\delta$ then consists in maximizing this smallest conic eigenvalue. Hence, we want to solve 
\begin{align}
\label{maxsigmin}
\max_{\delta \in \{0,1\}^L} \quad \sigma_{\min}(I_{\delta,L},
{\rm Tang}_{\mathcal C}(w_0)).
\end{align} 

One important remark to make at this point is that our criterion takes into account the underlying geometry of the estimation problem via the introduction of the 
tangent cone to $\mathcal C$ at $w_0$. 

\subsection{Preliminary experiments with a relaxation}
We performed some preliminary experiments with a small network. We relaxed the binary constraint to the convex set $[0,1]^L$ and used a level bundle method to optimise the relaxed problem

\begin{align}
\max_{\delta \in [0,1]^L} \quad \sigma_{\min}(I_{\delta,L},
{\rm Tang}_{\mathcal C}(w_0)).
\end{align} 

In the application, $w_0$ is clearly unknown and it has to be replaced with an appropriate estimator $\hat w_0$. For this purpose, we solved the estimation problem using different estimated values obtained from pseudo-measurements \cite{clements2011impact}. 

An example of the type of result we can obtain is given in the Figure \ref{resgrid} below. This example was obtained using a real network with 150 buses.

\begin{figure}[!th]
\includegraphics[width = 12cm]{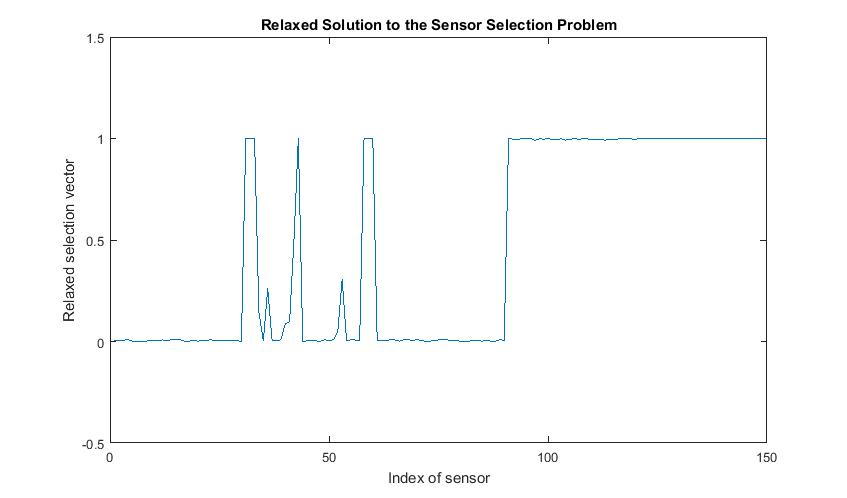}
\caption{\label{resgrid} The indicator vector for the sensor selection problem in a real world power grid network with 150 buses.} 
\end{figure}

This experiment shows that the relaxation based on a conic singular value is able to provide a simple method for sensor placement  that takes into account the constraints in the E-optimal-type of design. 
  
Thorough experiments with larger networks will be undertaken in future works on power grid, state estimation and sensor placement based on the smallest singular value criterion.

\section{Conclusion}

The goal of the present paper was to propose a fast algorithm for computing the smallest conic singular value of a matrix. The notion of conic singular value has appeared in recent works on compressed sensing as a very natural object to study and might be very useful in other applications as well. Our paper presents the first analysis of the numerical aspects of the smallest conic singular value. We found that the BFGS quasi-Newton method is running fast on simulated data. We will further experiment with real datasets from statistics, design of quadratic experiments and Compressed Sensing in future works.

	\bibliographystyle{plain}
	\bibliography{conic}
\end{document}